\documentclass[10pt]{amsart}
\usepackage{amssymb}
\usepackage{bbm} 
\usepackage{bm} 
\usepackage{graphicx}
\newcommand{\mc}{\mathcal}

\newcommand{\bb}{\mathbb}
\newcommand{\stackbounds}[2]{\stackrel{\mbox{\scriptsize #1}}{\mbox{\scriptsize #2}}}
\DeclareMathOperator{\im}{{\rm im}}

\DeclareMathOperator{\codim}{{\rm codim}}
\newcommand{\id}{{\rm id}}

\DeclareMathOperator{\spann}{{\rm span}}

\DeclareMathOperator{\inn}{{\rm in}}

\DeclareMathOperator{\vol}{{\rm vol}}
\DeclareMathOperator{\conv}{{\rm conv}}

\DeclareMathOperator{\face}{{\rm face}}
\DeclareMathOperator{\normal}{{\rm normal}}
\DeclareMathOperator{\vertex}{{\rm vertex}}
\DeclareMathOperator{\relint}{{\rm relint}}
\DeclareMathOperator{\lin}{{\rm lin}}
\DeclareMathOperator{\Trop}{{\rm Trop}}
\newcommand{\fan}{{\rm fan}} 
\newcommand{\unbal}{{\rm unbal}}
\newcommand{\refl}{{\rm refl}}
\DeclareMathOperator{\Poly}{{\rm Poly}}
\DeclareMathOperator{\Chow}{{\rm Chow}}
\DeclareMathOperator{\ChowSub}{{\rm Chow}_\nu} 
\DeclareMathOperator{\Chowi}{{\rm Chow}_{\cl\iota}}
\DeclareMathOperator{\ChowiSub}{{\rm Chow}_{\cl\iota,\nu}} 
\newcommand{\Gr}{\mathrm{Gr}}
\newcommand{\Hilb}{\mathrm{Hilb}}
\newcommand{\bbK}{\bb K} 
\newcommand{\bbk}{\mathbbm{k}} 
\newcommand{\Z}{\bb Z}

\newcommand{\R}{\bb R}
\newcommand{\Rp}{\bb R_{\geq 0}}

\newcommand{\lB}{\mathopen{\{\!\!\{}}
\newcommand{\rB}{\mathclose{\}\!\!\}}}
\newcommand{\msc}{\boxplus}
\newcommand{\Berg}{{\mc B}}
\newcommand{\coBerg}{{\mc B}^\ast}
\newcommand{\cl}{\bm} 
\newcommand{\eM}[1]{e^{#1}}
\newcommand{\eN}[1]{e_{#1}}

\DeclareRobustCommand{\qedify}[1]{%
  \ifmmode \quad\hbox{#1}
  \else
    \leavevmode\unskip\penalty9999 \hbox{}\nobreak\hfill
    \quad\hbox{#1}%
  \fi
}
\newtheorem{theorem}{Theorem}[section]
\newtheorem{maintheorem}[theorem]{Main theorem}
\newtheorem{lemma}[theorem]{Lemma}
\newtheorem{proposition}[theorem]{Proposition}
\newtheorem{corollary}[theorem]{Corollary}
\newtheorem{conjecture}[theorem]{Conjecture}
\theoremstyle{definition}
\newtheorem{definition}[theorem]{Definition}
\newtheorem{case*}[theorem]{Case}
\newenvironment{case}[1][]{\begin{case*}[#1]\pushQED{\qedify{$\diamondsuit$}}}{\popQED\end{case*}}
\newtheorem{example*}[theorem]{Example}
\newenvironment{example}[1][]{\begin{example*}[#1]\pushQED{\qedify{$\diamondsuit$}}}{\popQED\end{example*}}

\theoremstyle{remark}
\newtheorem{remark}[theorem]{Remark}

\numberwithin{equation}{section}

\begin{document}

\title{Tropical cycles and Chow polytopes} 
\author{Alex Fink}
\begin{abstract}
The Chow polytope of an algebraic cycle in a torus depends only on its tropicalisation.
Generalising this, we associate a Chow polytope to any abstract tropical variety
in a tropicalised toric variety.
Several significant polyhedra associated to tropical varieties are special cases of our Chow polytope.
The Chow polytope of a tropical variety $X$ is given by a simple combinatorial construction:
its normal subdivision is the Minkowski sum of $X$ and a reflected skeleton of
the fan of the ambient toric variety.
\end{abstract}
\maketitle

\section{Introduction}\label{sec:introduction}
Several well understood classes of tropical variety
are known to correspond to certain regular subdivisions of polytopes,
in a way that provides a bijection of combinatorial types.  

\begin{enumerate}
\item {\em Hypersurfaces} in~$\bb P^{n-1}$ 
are set-theoretically cut out by principal prime ideals.
If the base field has trivial valuation, then 
$\Trop\cl V(f)$ is\footnote{Throughout this paper we use boldface 
for classical algebro-geometric objects (except those with standard symbols
in blackboard bold or roman, which we preserve), and plain italic for tropical ones.} 
the fan of all cones of positive codimension
in the normal fan to its Newton polytope~${\rm Newt}(f)$. 
In the case of general valuation, the valuations of coefficients 
in $f$ induce a regular subdivision of~${\rm Newt}(f)$, and
$\Trop\cl V(f)$ consists of the non-full-dimensional faces in the
normal complex (in the sense of Section~\ref{sec:normal complex}).  

\smallskip
\item {\em Linear spaces} in~$\bb P^{n-1} = \bb P(\bbK^n)$ are cut out by ideals 
generated by linear forms.  
To a linear space $\cl X$ of dimension $n-d-1$ is associated a matroid $M(\cl X)$, whose bases
are the sets~$I\in\binom{[n]}d$ such that the projection of 
$\cl X$ to the coordinate subspace $\bbK\{e_i : i\not\in I\}$ has full rank.  
If the base field has trivial valuation, then $\Trop\bf X$ is a subfan 
(the {\em Bergman fan}~\cite{AK})
of the normal fan to the matroid polytope 
\begin{equation}\label{eq:matroid polytope}
\Poly(M(\cl X))=\conv\{\textstyle\sum_{j\in J}e_j : \mbox{$J$ is a basis of $M(\cl X)$}\} 
\end{equation}
of~$M(\cl X)$.
In the case of general valuations, the
valuations of the Pl\"ucker coordinates induce a regular subdivision 
of~$\Poly(M(\cl X))$ into matroid polytopes, and $\Trop\cl X$ 
consists of appropriate faces of the normal complex. 

\smallskip
\item {\em Zero-dimensional tropical varieties} are simply point configurations.
A zero-dimensional tropical variety $X$ 
is associated to an arrangement $\mc H$ of {\em upside-down} tropical hyperplanes
with cone points at the points of $X$: for instance, 
the tropical convex hull of the points of~$X$ is a union of closed regions
in the polyhedral complex determined by $\mc H$.  
The arrangement $\mc H$ is dual to a fine mixed subdivision of a simplex,
and $X$ consists of the faces dual to little simplices in the normal complex
of this subdivision.
\end{enumerate}

The polytopes and subdivisions in this list are special cases of the {\em Chow polytope},
or subdivision of Chow polytopes,
associated to any cycle $\cl X$ on~$\bb P^{n-1}$ as 
the weight polytope of the point representing $\cl X$ in the {\em Chow variety},
the parameter space of cycles.
Although this is an entirely classical construction, in fact
the Chow polytope subdivision of $\cl X$ depends only on the tropical variety $\Trop\cl X$, 
and the construction can be extended to associate Chow polytope
subdivisions to all tropical varieties in~$\bb R^{n-1}$.

This paper's main theorem, Theorem~\ref{r:1}, provides
a simple tropical formula for this Chow polytope subdivision in
terms of $\Trop\cl X$, making use of a {\em stable Minkowski sum} operation on tropical cycles
introduced in Section~\ref{sec:msc}.  The formula is similar to its
classical analogue, and is even simpler in one salient respect, namely that there's
no need to invoke any sort of Grassmannian (Remark~\ref{rem:cl vs trop}).
The formula generalises to subvarieties of any projective toric variety.

There is however no general map in the reverse direction, from Chow polytope subdivision
to tropical variety (that is, the bijection of combinatorial types in the opening
examples is a special phenomenon).  In Section~\ref{sec:kernel} we present an example of
two distinct tropical varieties with the same Chow polytope.

Finally, in Section~\ref{sec:linear} we use this machinery to at last
record a proof of the fact that tropical linear spaces are exactly tropical varieties of degree~1.

\section{Tropical setup}

We begin with a few polyhedral notations and conventions.
For $\Pi$ a polyhedron in a real vector space~$V$ and $u:V\to\R$ 
a linear functional, $\face_u\Pi$ is the face of~$\Pi$ on which
$u$ is minimised, if such a face exists.  
{\newcommand{\Rho}{\mathrm P}
For $\Pi,\Rho$ polyhedra, $\Pi+\Rho$ is the Minkowski sum
$\{\pi+\rho:\pi\in\Pi, \rho\in\Rho\}$, and we write
$-\Rho=\{-\rho:\rho\in\Rho\}$ and $\Pi-\Rho = \Pi+(-\Rho)$.
}

\subsection{Tropical cycles}\label{ssec:tropical cycles}
Let $N_\R$ be a real vector space 
containing a distinguished full-dimensional lattice $N$, so that
$N_\R = N\otimes\R$.  This is all the structure necessary to define abstract tropical
cycles in $N_\R$, and this is the context in which we will work at first.
However, we will often have the situation of Case~\ref{case:projective}.

\begin{case}[Projective tropical varieties]\label{case:projective}
Let $\cl X$ be a classical subvariety of $\bb P^{n-1}$ tropicalised with respect to the torus 
$(\bbK^\ast)^n/\bbK^\ast\subseteq\bb P^{n-1}$, where the $\bbK^\ast$ in the quotient embeds diagonally.
Then $X:=\Trop\cl X$ is a tropical fan in~$N_\R=\bb R^n/(1,\dots,1)$,
and $N=\bb Z^n/(1,\ldots,1)$ is the lattice of integer points within~$N_\R$.
The dual vector space to $N_\R$ is $M_\R=(1,\ldots,1)^\perp=(N_\R)^\vee$
(sometimes it will be convenient to use a translate instead).
This $M_\R$ also carries its lattice $M = M_\R\cap\bb Z^n = N^\vee$.  

For maximal clarity we will write $\eN{i}$ for the image in $N_\R$ of
a basis element of~$\bb R^n$, and $\eM{i}$ for a basis element in the $(\bb R^n)^\ast$ 
of which $M_\R$ is a subspace.
For $J\subseteq[n]$, the notation $\eN{J}$ means $\sum_{j\in J}\eN{j}$, and
$\eM{J}$ is analogously defined.
\end{case}

For a polyhedron $\sigma\subseteq N_\R$,
let $\lin\sigma$ be the translate of the affine hull of $\sigma$ to the origin.
We say that $\sigma$ is {\em rational} if $N_\sigma := N\cap\lin\sigma$ is
a lattice of rank~$\dim\sigma$.

The fundamental tropical objects we will be concerned with 
are abstract {\em tropical cycles} in $N_\R$.  
See \cite[Section~5]{AR} for a careful exposition
of tropical cycles.  Loosely, a tropical cycle $X$ of dimension~$k$ consists of
the data of a rational polyhedral complex $\Sigma$ pure of dimension~$k$, 
and for each facet $\sigma$ of~$\Sigma$ an integer multiplicity $m_\sigma$ satisfying 
a {\em balancing condition} at codimension~1 faces,
modulo identifications which ensure that the precise choice 
of polyhedral complex structure, among those with a given support, is unimportant.  
A {\em tropical variety} is an effective tropical cycle, one in which all multiplicities
$m_\sigma$ are nonnegative.

We write $ Z_k$ for the additive group of tropical cycles in $N_\R$ of dimension~$k$.
We also write $ Z= \bigoplus_k Z_k$, and use upper indices for codimension,
$ Z^k = Z_{\dim N_\R-k}$.
If $\Sigma$ is a polyhedral complex,
then by $ Z(\Sigma)$ (and variants with superscript or subscript)
we denote the group of tropical cycles $X$ (of appropriate dimension) which can be given
some polyhedral complex structure with underlying polyhedral complex $\Sigma$.
Our notations $Z$ and $Z^k$ are compatible with~\cite{AR},
but we use $Z(\Sigma)$ differently (in~\cite{AR} it refers merely to cycles contained as sets
in~$\Sigma$, a weaker condition).

If a tropical cycle $X$ can be given a polyhedral complex structure which is a 
fan over the origin, we call it a {\em fan cycle}.
We prefer this word ``fan'', as essentially in~\cite{GKM},
over ``constant-coefficient'', for brevity and for not suggesting tropicalisation;
and over the ``affine'' of~\cite{AR}, since {\em tropical affine space}
should refer to a particular partial compactification of~$N_\R$.
We use notations based on the symbol $ Z^\fan$ for groups of tropical fan cycles.

In a few instances it will be technically convenient to work with objects which 
are like tropical cycles except that the balancing condition is not required.
We call these {\em unbalanced cycles} and use notations based on the symbol $Z^\unbal$.  
That is, $Z^\unbal$ simply denotes the free Abelian group on the cones of~$\Delta$.
If $\sigma\subseteq N_\R$ is a $k$-dimensional polyhedron, we write
$[\sigma]$ for the unbalanced cycle $\sigma$ bearing multiplicity~1,
and observe the convention $[\emptyset]=0$.  
Then every tropical cycle can be written as an integer combination of various~$[\sigma]$.

It is a central fact of tropical intersection theory that $ Z^\fan$ is a graded ring,
with multiplication given by ({\em stable}) {\em tropical intersection},
which we introduce next, and grading given by codimension.  
The invocation of these notions in the toric context~\cite[Section~4]{FS}
prefigured certain aspects of the tropical machinery:
\begin{theorem}[Fulton--Sturmfels]
Given a complete fan $\Sigma$, $ Z^\fan(\Sigma)$ is the Chow cohomology ring
of the toric variety associated to $\Sigma$.
\end{theorem}

Given two rational polyhedra $\sigma$ and~$\tau$,
we define a multiplicity $\mu_{\sigma,\tau}$ 
arising from the lattice geometry, namely the index
$$\mu_{\sigma,\tau} = [N_{\sigma+\tau} : N_\sigma + N_\tau].$$
We define two variations where we require, respectively, transverse intersection
and linear independence:
\begin{align*}
\mu^\bullet_{\sigma,\tau} &= \left\{\begin{array}{ll}
\mu_{\sigma,\tau} & \mbox{if }\codim(\sigma\cap\tau) = \codim\sigma + \codim\tau \\ 
0 & \mbox{otherwise,} 
\end{array}\right. \\
\mu^\msc_{\sigma,\tau} &= \left\{\begin{array}{ll}
\mu_{\sigma,\tau} & \mbox{if }\dim(\sigma+\tau) = \dim\sigma + \dim\tau \\ 
0 & \mbox{otherwise.} 
\end{array}\right.
\end{align*}

Alternatively, $\mu^\msc_{\sigma,\tau}$ is the absolute value of the
determinant of a block matrix consisting of
a block whose rows generate $N_\sigma$ as a $\Z$-module above a block whose rows generate~$N_\tau$,
in coordinates providing a basis for any $(\dim\sigma + \dim\tau)$-dimensional lattice containing
$N_{\sigma+\tau}$.  Likewise $\mu^\bullet_{\sigma,\tau}$ can be computed
from generating sets for the dual lattices.

If $\sigma$ and~$\tau$ are polytopes in~$N_\R$ which are either disjoint or
intersect transversely in the relative interior of each, their stable tropical intersection is
\begin{equation}\label{eq:trop int of facets}
[\sigma]\cdot[\tau] = \mu^\bullet_{\sigma,\tau}[\sigma\cap\tau].
\end{equation}
If  $X = \sum_\sigma m_\sigma[\sigma]$ 
and $Y = \sum_\tau   n_\tau  [\tau  ]$ 
are unbalanced cycles such that every pair of facets $\sigma$ of~$X$
and $\tau$ of~$Y$ satisfy this condition,
then their stable tropical intersection is obtained by linear extension,
\begin{equation}\label{eq:generic trop int}
X\cdot Y=\sum_{\sigma,\tau}
m_\sigma n_\tau\cdot \mu^\bullet_{\sigma,\tau}[\sigma\cap\tau].
\end{equation}
If $X$ and $Y$ are tropical cycles, so is~$X\cdot Y$ (see~\cite{AR}).
For a point $v\in N_\R$, let $[v]\msc Y$ denote the translation of~$Y$ by~$v$;
this is a special case of a notation we introduce in Section~\ref{sec:msc}.
If $X$ and~$Y$ are rational tropical cycles with no restrictions, 
then for generic small displacements $v\in N_\R$ the
faces of $X$ and~$[v]\msc Y$ intersect suitably for 
equation~\eqref{eq:generic trop int} to be applied.
In fact the facets of the intersection $X\cdot([v]\msc Y)$ 
vary continuously with~$v$, in a way that can be continuously extended to all $v$.
This is essentially the {\em fan displacement rule} of~\cite{FS}, which ensures that
$X\cdot Y$ is always well-defined.
\begin{definition}\label{def:trop int}
Given two tropical cycles $X,Y$, their ({\em stable}) {\em tropical intersection} is 
$$X\cdot Y=\lim_{v\to 0} X\cdot([v]\msc Y).$$  
\end{definition}

We introduce a few more operations on cycles.
Firstly, there is a cross product defined in the expected fashion.
Temporarily write $Z(V)$ for the ring of tropical cycles defined in the vector space $V$.
Let $(N_i)_\R$, $i=1,2$, be two real vector spaces.
Then there is a well-defined bilinear cross product map 
$$\times:Z^\unbal((N_1)_\R)\otimes Z^\unbal((N_2)_\R)\to Z^\unbal((N_1\oplus N_2)_\R)$$
linearly extending $[\sigma]\times[\tau]=[\sigma\times\tau]$, and the exterior product
of tropical cycles is a tropical cycle.

Let $h:N\to N'$ be a linear map of lattices, inducing a map of real vector spaces
which we will also denote $h:N_\R\to N'_\R$ (an elementary case of a tropical morphism).  
Cycles can be pushed forward and pulled back along~$h$.  These are special cases of notions
defined in tropical intersection theory even in ambient tropical varieties 
other than $\bb R^n$ (in the general case, one can push forward general cycles
but only pull back complete intersections of Cartier divisors~\cite{AR}).

Given a cycle
$Y=\sum_\sigma m_\sigma[\sigma]$ on~$N'_\R$, its pullback is
defined in~\cite{Allermann} as follows.  This is shown in~\cite[Proposition~2.7]{FS} to agree
with the pullback on Chow rings of toric varieties.
$$h^\ast(Y) = 
\hspace{-1em}\sum_{\sigma\ :\ \mbox{\scriptsize$\sigma$ meets $\im h$ transversely}}\hspace{-1em}
m_\sigma [N_{h^{-1}(\sigma)} : h^{-1}(N'_\sigma)] [h^{-1}(\sigma)]$$
The pushforward is defined in~\cite{GKM} in the tropical context,
and is shown to coincide with the cohomological pushforward in~\cite[Lemma~4.1]{Katz}.
If $X=\sum_\sigma m_\sigma[\sigma]$ is a cycle on~$N_\R$, 
its pushforward is
$$h_\ast(X) = \sum_{\sigma\ :\ \mbox{\scriptsize$h|_\sigma$ injective}}
m_\sigma [N'_{h(\sigma)} : h(N_\sigma)] [h(\sigma)].$$
In these two displays, the conditions on~$\sigma$ in the sum are equivalent to
$h^{-1}(\sigma)$ or $h(\sigma)$, respectively, having the expected dimension.
Pushforwards 
and pullbacks of tropical cycles are
tropical cycles.  

\subsection{Normal complexes}\label{sec:normal complex}

Write $M = N^\vee$, $M_\R = N_\R{}^\vee$ for the dual lattice and real vector space.
Let $\pi:M_\R\times\R\to M_\R$ be the projection to the first factor.  
A polytope $\Pi\subseteq\bb M_\R\times\R$
induces a regular subdivision $\Sigma$ of $\pi(\Pi)$.
Our convention will be that regular subdivisions are determined by lower
faces: so the faces of~$\Sigma$ are the projections $\pi(\face_{(u,1)}\Pi)$. 
We will also write $\face_u \Sigma$ to refer to this last face.  
In general, we will not consider regular subdivisions $\Sigma$ 
by themselves but will also want to retain the data of~$\Pi$.  
More precisely, what is necessary is to have a well-defined normal complex;
for this we need only $\Sigma$ together with 
the data of the heights of the vertices of~$\Pi$ visible from underneath, 
equivalently the lower faces of~$\Pi$.  (When we refer to ``vertex heights'' we shall
always mean only the lower vertices.)

\begin{definition} 
The ({\em inner}) {\em normal complex} $\mc N(\Sigma,\Pi)$ to
the regular subdivision $\Sigma$ induced by~$\Pi$ is
the polyhedral subdivision of~$N_\R$ with a face 
$$\normal(F) = \{u\in N_\R : W\subseteq \face_{(u,1)}(\Pi)\}$$
for each face $F=\conv(\pi(W))$ of~$\Sigma$.
\end{definition}
We will allow ourselves to write $\mc N^1(\Sigma)$ for $\mc N^1(\Sigma,\Pi)$
when $\Pi$ is clear from context.  
If $\Pi$ is contained in $M_\R\times\{0\}$, which we identify with~$M_\R$,
then $\mc N(\Sigma,\Pi)$ is the {\em normal fan} of~$\Pi$.

We give multiplicities to the faces of the skeleton $\mc N^e(\Sigma,\Pi)$
of $\mc N(\Sigma,\Pi)$ so as to make it a cycle, 
which we also denote $\mc N^e(\Sigma,\Pi)$.
To each face $\normal(F)\in\mc N(\Sigma,\Pi)$ of codimension~$e$,
we associate the multiplicity $m_{\normal(F)} = \vol F$ where $\vol$ is the 
normalised lattice volume, i.e. the Euclidean volume
on~$\lin F$ rescaled so that
any simplex whose edges incident to one vertex form a basis for $N_F$ has volume~1.
In fact $\mc N^e(\Sigma,\Pi)$ is a tropical cycle.
In codimension~1 a converse holds as well.

\begin{theorem}\label{r:weights on polytopes}\mbox{}
\begin{enumerate}
\renewcommand{\labelenumi}{(\alph{enumi})}
\item For any rational regular subdivision $\Sigma$ in $M_\R$
induced by a polytope $\Pi$ in $M_\R\times\R$, the skeleton
$\mc N^e(\Sigma,\Pi)$ is a tropical variety.
\item For any tropical variety $X\in Z^1(N_\R)$, 
there exists a rational polytope $\Pi$ in $M_\R\times\R$ and induced regular subdivision $\Sigma$,
unique up to translation and adding a constant to the vertex heights, 
such that $X=\mc N^1(\Sigma,\Pi)$.
\end{enumerate}
\end{theorem}
Part~(a) in the case of fans, i.e.\ $\Pi\subseteq M_\R\times\{0\}$, 
is a foundational result in the polyhedral 
algebra~\cite[Section~11]{McMullen2}.
The statement for general tropical varieties follows since the 
normal complex of~$\Sigma$ is just the slice through the normal fan of~$\Pi$ at height~1,
and this slicing preserves the balancing condition.
Part~(b) is also standard, and is a consequence of
ray-shooting algorithms, the codimension~1 case of~Theorem~\ref{r:DFS2.2}.

One more fact will be important when we move beyond $\bb P^{n-1}$ as ambient variety.
This is the content of~\cite[Theorem 5.1]{McMullen2} cast tropically. 
\begin{theorem}\label{r:slice-project}
Let $\iota:N\to N'$ be an inclusion of lattices such that $\iota N$ is saturated in~$N'$,
and $\iota^{\mathrm T}$ the dual projection.  
For any polytope $\Pi'$ in $M'_\R\times\R$, let $\Pi=(\iota^{\mathrm T}\times\id)\Pi'$
be its projection to $M_\R\times\R$, and let $\Sigma'$ and $\Sigma$ be the induced regular subdivisions.
Then $\mc N^e(\Sigma) = \mc N^e(\Sigma')\cdot [\iota N]$.
\end{theorem}
This $\Sigma$ is the {\em image subdivision} of~$\Sigma'$ of~\cite{KatzToolkit}; 
this is the natural notion of projection for regular subdivisions with vertex heights.

\section{Minkowski sums of cycles}\label{sec:msc}

Let $N$ be any lattice.  
For a tropical cycle $X=\sum m_\sigma[\sigma]$, 
we let $X^\refl=\sum m_\sigma[-\sigma]$ denote its reflection about the origin.
(This is the pushforward or pullback of~$X$ along the linear isomorphism $x\mapsto -x$.)  

Given two polyhedra $\sigma,\tau\subseteq N_\R$, define the ({\em stable}) {\em Minkowski sum}
\begin{equation}\label{eq:def msc}
[\sigma]\msc[\tau] = \mu^\msc_{\sigma,\tau}[\sigma+\tau].
\end{equation}
%
Compare \eqref{eq:trop int of facets}.
If $X$ and~$Y$ are cycles in $N_\R$, then we can
write their intersection and Minkowski sum in 
terms of their exterior product $X\times Y\in N_\R\times N_\R$.  
We have an exact sequence
$$0\to N_\R\stackrel\iota\to N_\R\times N_\R\stackrel\phi\to N_\R\to 0$$
of vector spaces where $\iota$ is the inclusion along the diagonal
and $\phi$ is subtraction, $(x,y)\mapsto x-y$.  It is then routine
to check from the definitions that 
\begin{align}\label{eq:push and pull}
X\cdot Y &= \iota^\ast(X\times Y) \notag\\
X\msc Y^\refl &= \phi_\ast(X\times Y) 
\end{align}
Since pullback is well-defined and takes tropical cycles to tropical cycles,
it follows immediately that there is a well-defined bilinear map
${\msc}: Z^\unbal\otimes Z^\unbal\to Z^\unbal$ extending~\eqref{eq:def msc}, 
restricting to a bilinear map ${\msc}: Z\otimes  Z\to  Z$.

A notion of Minkowski sum for tropical varieties arose in~\cite{CMS}
as the tropicalisation of the {\em Hadamard product} for classical varieties.
The Minkowski sum of two tropical varieties in that paper's 
sense can have dimension less than the expected dimension.
By contrast our bilinear operation $\msc$ should be regarded as a {\em stable} Minkowski sum
for tropical cycles.  It is additive in dimension, i.e.\
$ Z_d\msc Z_{d'}\subseteq Z_{d+d'}$, just as 
stable tropical intersection is additive in codimension.  
The next lemma further relates intersection and Minkowski sum.

The balancing condition implies that for
any tropical cycle $X$ in~$N_\R$ of dimension~$\dim N_\R$, $X(u)$ is constant
for any $u\in N_\R$ for which it's defined.  We shall denote this constant $\deg X$.
Similarly, if $\dim X=0$, then $X$ is a finite sum of points with multiplicities, 
and we will let $\deg X$ be the sum of these multiplicities.  
These are both special cases of Definition~\ref{def:degree}, to come.

\begin{lemma}\label{r:two pairings}
Let $X$ and~$Y$ be tropical cycles on~$N_\R$, of complementary dimensions.  Then 
$$\deg(X\cdot Y) = \deg(X\msc Y^\refl).$$
\end{lemma}

\begin{proof}
Let $u\in N_\R$ be generic.  Let
$\Sigma(X)$ and~$\Sigma(Y)$ be polyhedral complex structures on $X$ and~$Y$.
The multiplicity of~$X\msc Y^\refl$ at a point~$u\in N_\R$ is
\begin{equation}\label{eq:XY}
(X\msc Y^\refl)(u) = \sum_{\sigma,\tau} \mu^\msc_{\sigma,\tau},
\end{equation}
summing over only those $\sigma\in\Sigma(X)$ and
$\tau\in\Sigma(Y)^\refl$ with
$u\in\sigma+\tau$, i.e. with
$(\{u\}-\tau)\cap\sigma$ nonempty.  
These $\{u\}-\tau$ are the cones of $\Sigma(Y')$, where $Y' = [u]\msc Y$.
Then by~\eqref{eq:generic trop int}, $\deg(X\cdot Y')$
is given by the very same expression~\eqref{eq:XY} except with $\mu^\bullet$ in place of~$\mu^\msc$;
and by the fan displacement rule preceding Definition~\ref{def:trop int}, 
$\deg(X\cdot Y)=\deg(X\cdot Y')$.
But $\mu^\bullet_{\sigma,\tau}=\mu^\msc_{\sigma,\tau}$ when $\sigma$ and~$\tau$
are of complementary dimensions.
\end{proof}

\begin{lemma}\label{r:slicing commutes}
Let $X$ and $Z$ be tropical cycles on~$N_\R$, and $Y$ a cycle which is a classical linear space 
through the origin, with $X\subseteq Y$.  Then
$$X\msc(Y\cdot Z) = Y\cdot(X\msc Z).$$
\end{lemma}

\begin{proof}
Replacing $Z$ (and thus $X\msc Z$) by a generic small translate, we may take the
intersections to be set-theoretic intersections with lattice multiplicity. 
By linearity, we may assume $X$ and $Z$ are of the form $[\sigma]$.
Then this reduces to checking set-theoretic equality and checking equality
of multiplicities, both of which are routine.
\end{proof}

We specialise to Case~\ref{case:projective}.
Let $\mc L$ be the fan
of the ambient toric variety~$\bb P^n$, which is the normal fan in~$N$
to the standard simplex $\conv\{\eM{i}\}$.  The ray generators
of~$\mc L$ are $\eN{i}\in N$, and every proper subset of the rays
span a face, which is simplicial.  For $J\subsetneq[n]$ 
let $C_J=\Rp\{\eN{j}:j\in J\}$ be the face of~$\mc L$
indexed by~$J$.
Let $\mc L_k$ be the dimension~$k$ skeleton of~$\mc L$ with multiplicities 1,
that is, the canonical $k$-dimensional tropical fan linear space.  

\begin{definition}[{\cite[Definition~9.13]{AR}}]\label{def:degree}
The {\em degree} of a tropical cycle $X\in  Z^e(N_\R)$ is
$\deg X := \deg(X\cdot\mc L_e)$.
\end{definition}

The symbol $\deg$ appearing on the right side is the special case defined 
just above for cycles of dimension 0.
It is a consequence of the fan displacement rule that $\deg X = \deg(X\cdot([v]\msc\mc L_e))$
for any $v\in N_\R$.

\begin{lemma}\label{r:ch deg}
Let $X\in  Z^e$.  Then
$$\deg(X\msc{\mc L_{e-1}}^\refl) = e\deg X.$$
\end{lemma}

\begin{proof}
By Lemma~\ref{r:two pairings} we have
\begin{align*}
\deg(X\msc{\mc L_{e-1}}^\refl)
  &= \deg((X\msc{\mc L_{e-1}}^\refl)\cdot\mc L_1)
\\&= \deg(X\msc{\mc L_{e-1}}^\refl\msc{\mc L_1}^\refl)
\\&= \deg(({\mc L_{e-1}}^\refl\msc{\mc L_1}^\refl)\cdot X^\refl)
\\&= \deg((\mc L_{e-1}\msc\mc L_1)\cdot X)
\\&= \deg((e\mc L_e)\cdot X)
\\&= e\deg(X\cdot\mc L_e)
\\&= e\deg X.\qedhere
\end{align*}
\end{proof}

\begin{remark}
The classical projection formula of intersection theory is valid tropically \cite[Proposition 7.7]{AR},
and has an analogue for~$\msc$.
For a linear map of lattices $h:N\to N'$ and cycles 
$X\in Z(N_\R)$ and $Y\in Z(N'_\R)$, 
we have
\begin{align*}
h_\ast(X\cdot h^\ast(Y)) &= h_\ast(X)\cdot Y,\\
X\msc h^\ast(Y) &= h^\ast(h_\ast(X)\msc Y).
\end{align*}

The facts in this section, as well as 
the duality given by polarisation in the algebra of cones which 
exchanges intersection and Minkowski sum,
are all suggestive of the existence of a duality between tropical stable 
intersection and stable Minkowski sum.
However, we have not uncovered a better statement of such a duality than 
equations~\eqref{eq:push and pull}.
\end{remark}

\section{Chow polytopes}\label{sec:Chow}
In this section we introduce Chow polytopes.  There is little new content here:
see \cite{KSZ}, \cite[ch.~4]{GKZ} and~\cite{Dalbec} for fuller treatments of 
this material, the first for the toric background, the second 
in the context of elimination theory, and 
the last especially from a computational standpoint.
The assumptions of Case~\ref{case:projective} will be in force for most of this section,
and most of the rest of the paper.

Let $\bbK$ be an algebraically closed field. 
Let $(\bbK^\ast)^n$ be an algebraic torus acting via a linear representation on a vector space~$V$, 
or equivalently on its projectivisation $\bb P(V)$.  
Suppose that the action of $(\bbK^\ast)^n$ is {\em diagonalisable},
i.e.\ $V$ can be decomposed as a direct sum $V=\bigoplus V_i$ where 
$(\bbK^\ast)^n$ acts on each $V_i$ by a character or {\em weight} $\chi^{w_i}:(\bbK^\ast)^n\to\bb K^\ast$.
A character $\chi^{w_i}$
corresponds to a point $w_i$ in the character lattice of~$(\bbK^\ast)^n$,
via $\chi^{w_i}(t) = t^{w_i}$.
We shall always assume $V$ is finite-dimensional, except in a few instances
where we explicitly waive this assumption for technical convenience.   
If $V$ is finite-dimensional, the action of~$(\bbK^\ast)^n$ is necessarily diagonalisable.
\begin{definition}
Given a point $v\in V$
of the form $v=\sum_{k\in K} v_{i_k}$ with each $v_{i_k}\in V_{i_k}$ nonzero,
the {\em weight polytope} of~$v$ is $\conv\{w_k : k\in K\}.$
\end{definition}
\noindent 
If $\cl X\subseteq\bb P(V)$ is a $(\bbK^\ast)^n$-equivariant subvariety, this defines 
the weight polytope of a point $x\in\cl X$. 

The {\em Chow variety} $\Gr(d,n,r)$ of~$\bb P^{n-1}$, 
introduced by Chow and van der Waerden in~1937 \cite{CvdW},
is the parameter space of effective cycles of dimension~$d-1$
and degree~$r$ in $\bb P^{n-1}$.
When we invoke homogeneous coordinates on $\bb P^{n-1}$ 
we will name them $x_1,\ldots,x_n$.  
\begin{example}\mbox{}
\begin{enumerate}
\item The variety $\Gr(n-1,n,r)$ parametrising degree $r$ cycles of codimension 1
is $\bb P(\bbK[x_1,\ldots,x_n]_r)\cong\bb P^{\binom{r+n-1}r-1}$.
An irreducible cycle is represented by its defining polynomial.
\item The variety $\Gr(d,n,1)$ parametrises degree~1 effective cycles,
which must be irreducible and are therefore linear spaces.
So $\Gr(d,n,1)$ is simply the Grassmannian $\Gr(d,n)$,
motivating the notation.\qedhere
\end{enumerate}
\end{example}

The Chow variety $\Gr(d,n,r)$ is projective.  Indeed,
we can present the coordinate ring of $\Gr(n-d,n)$
in terms of (primal) {\em Pl\"ucker coordinates}, which we write as brackets:
$$\bbK[\Gr(n-d,n)] = \bbK\Bigl[\,[J] : J\in{\binom{[n]}{n-d}}\,\Bigr] 
\big/ (\mbox{Pl\"ucker relations}).$$
For our purposes the precise form of the Pl\"ucker relations will be unimportant.
Then $\Gr(d,n,r)$ has a classical embedding into 
the space $\bb P(\bbK[\Gr(n-d,n)]_r)$ of homogeneous degree~$r$
polynomials on~$\Gr(n-d,n)$ up to scalars, given by the {\em Chow form}~\cite{CvdW}.
We denote the Chow form of~$\cl X$ by~$R_{\cl X}$.

\begin{remark}\label{rem:meet}
For $\cl X$ irreducible, the Chow form $R_{\cl X}$ is the defining polynomial of
the locus of linear subspaces of~$\bb P^{n-1}$ of dimension~$n-d-1$ which intersect~$\cl X$.
There is a single defining polynomial since $\mathrm{Pic}(\Gr(n-d,n))=\bb Z$.
\end{remark}

The natural componentwise action
$(\bbK^\ast)^n\curvearrowright\bbK^n$ induces an action
$(\bbK^\ast)^n\curvearrowright S^\ast(\bigwedge^{n-d}\bbK^n)$.
The ring $\bbK[\Gr(n-d,n)]$ is a quotient of this symmetric algebra
by the ideal of Pl\"ucker relations.  This ideal is homogeneous in the weight grading, 
so the quotient inherits an $(\bbK^\ast)^n$-action.
The Chow variety is an $(\bbK^\ast)^n$-equivariant subvariety of~$\bb P(\bbK[\Gr(n-d,n)]_r)$,
so we also get an action $(\bbK^\ast)^n\curvearrowright\Gr(n-d,n)$.
The weight spaces of $\bbK[\Gr(n-d,n)]$ under the $(\bbK^\ast)^n$-action
are spanned by monomials in the brackets $[J]$.  
The weight of a bracket monomial $\prod_i [J_i]^{m_i}$ is $\sum_i m_i \chi^{J_i}$.

\begin{definition}
If $\cl X$ is a cycle on~$\bb P^{n-1}$
represented by the point $x$ of~$\Gr(d,n,r)$, the {\em Chow polytope} $\Chow(\cl X)$ of~$\cl X$ is 
the weight polytope of~$x$.
\end{definition}

\pagebreak[3] 
\begin{example}\label{ex:motivating2}\mbox{}
\begin{enumerate}
\item The Chow form of a hypersurface $\cl V(f)$ is simply its defining polynomial $f$
with the variables $x_k$ replaced by brackets $[k]$,
so that the Chow polytope $\Chow(\cl V(f))$ is the Newton polytope~${\rm Newt}(f)$.
\smallskip 
\item The Chow form of a $(d-1)$-dimensional linear space $\cl X$ is a linear form in the brackets,
$\sum_J p_J[J]$, where the $p_J$ are the {\em dual} Pl\"ucker coordinates of~$\cl X$
for $J\in\binom{[n]}{n-d}$.  
Accordingly $\Chow(\cl X)$ is the polytope $\Poly(M(\cl X)^\ast)$ of the dual matroid.
Note that this is simply the image of $\Poly(M(\cl X))$ under a reflection. 
\smallskip
\item For $\cl X=\cl X_A$ an embedded {\em toric variety} in $\bb P^{n-1}$, 
the Chow polytope $\Chow(\cl X)$ is the {\em secondary polytope} associated to the
vector configuration~$A$ \cite[Chapter~8.3]{GKZ}.\qedhere
\end{enumerate}
\end{example}

From a tropical perspective, the preceding setup has all pertained to the
constant-coefficient case.
Suppose now that the field $\bbK$ has a nontrivial valuation $\nu:\bbK^\ast\to\bb Q$,
with residue field $\bbk\hookrightarrow\bbK$.  
For instance we might take $\bbK=\bbk\lB t\rB$ 
the field of Puiseux series over an algebraically closed field~$\bbk$, with the valuation
$\nu:\bbK^\ast\to\bb Q$ by least degree of~$t$.  
Let $\cl X$ be a cycle on~$\bb P^{n-1}$ with
Chow form~$R_{\cl X}\in\bbK[\Gr(n-d,n)]$.
Let $\tau_1,\ldots,\tau_m\in\bbK$ be the coefficients of $R[X]$,
so that $R[X]$ is defined over the subfield $\bbk[\tau_1^{\pm1},\ldots,\tau_n^{\pm 1}]\subseteq\bbK$.
The restriction of $\nu$ to this subfield is a discrete valuation,
so we may assume that all the $\nu(\tau_i)$ are integers.

The torus $(\bbk^\ast)^n$ acts on~$\bbk[\Gr(n-d,n)]$ just as before,
and therefore acts on 
$\bbk[\Gr(n-d,n)][\tau_1^{\pm1},\ldots,\tau_n^{\pm 1}]$.  
Let $(\bbk^\ast)^n\times\bbk^\ast\curvearrowright \bbk[\tau^{\pm1}][\Gr(n-d,n)]$ 
where the right factor acts on Laurent monomials in $\tau_1,\ldots,\tau_n$, 
with $\tau^a$ having weight $\sum_{i=1}^m a_i\nu(\tau_i)$.  
Let $\Pi$ be the weight polytope of the Chow form~$R_{\cl X}$ with respect to this action.
\begin{definition}
The {\em Chow subdivision} of a cycle $\cl X$ on~$\bb P^{n-1}$
over $(\bbK,\nu)$ is the regular subdivision $\ChowSub(\cl X)$ induced by~$\Pi$.
\end{definition}

The Chow subdivision
is the non-constant-coefficient analogue of the Chow polytope, generalising the 
polytope subdivision of the opening examples.
It appears as the {\em secondary subdivision} in Definition~5.5 of~\cite{KatzToolkit}, 
but nothing is done with the definition in that work, and we believe this paper
is the first study to investigate it in any detail.
Observe that $\ChowSub(\cl X)$ is a subdivision of $\Chow(\cl X)$, and
if $\nu$ is the trivial valuation, $\ChowSub(\cl X)$ is $\Chow(\cl X)$ unsubdivided. 
By $\mc N(\ChowSub(\cl X))$ we will always mean $\mc N(\ChowSub(\cl X),\Pi)$.

If $(u,v):\bbk^\ast\to (\bbk^\ast)^n\times\bbk^\ast$ is a one-parameter subgroup which
as an element of $N\times\Z$ has negative last coordinate, 
then $\face_u \ChowSub(\cl X)=\face_{(u,v)} \Pi$ is bounded.  
We observe that a bounded face $F=\face_u \ChowSub(\cl X)$ of~$\ChowSub(\cl X)$ is 
the weight polytope of the toric degeneration $\lim_{t\to 0}u(t)\cdot\cl X$.
This follows from an unbounded generalisation of 
Proposition~1.3 of~\cite{KSZ}, 
which describes the toric degenerations of a point in terms of
the faces of its weight polytope.

\begin{example}\label{ex:Chow subdivision}
\begin{figure}
\centering{
\includegraphics{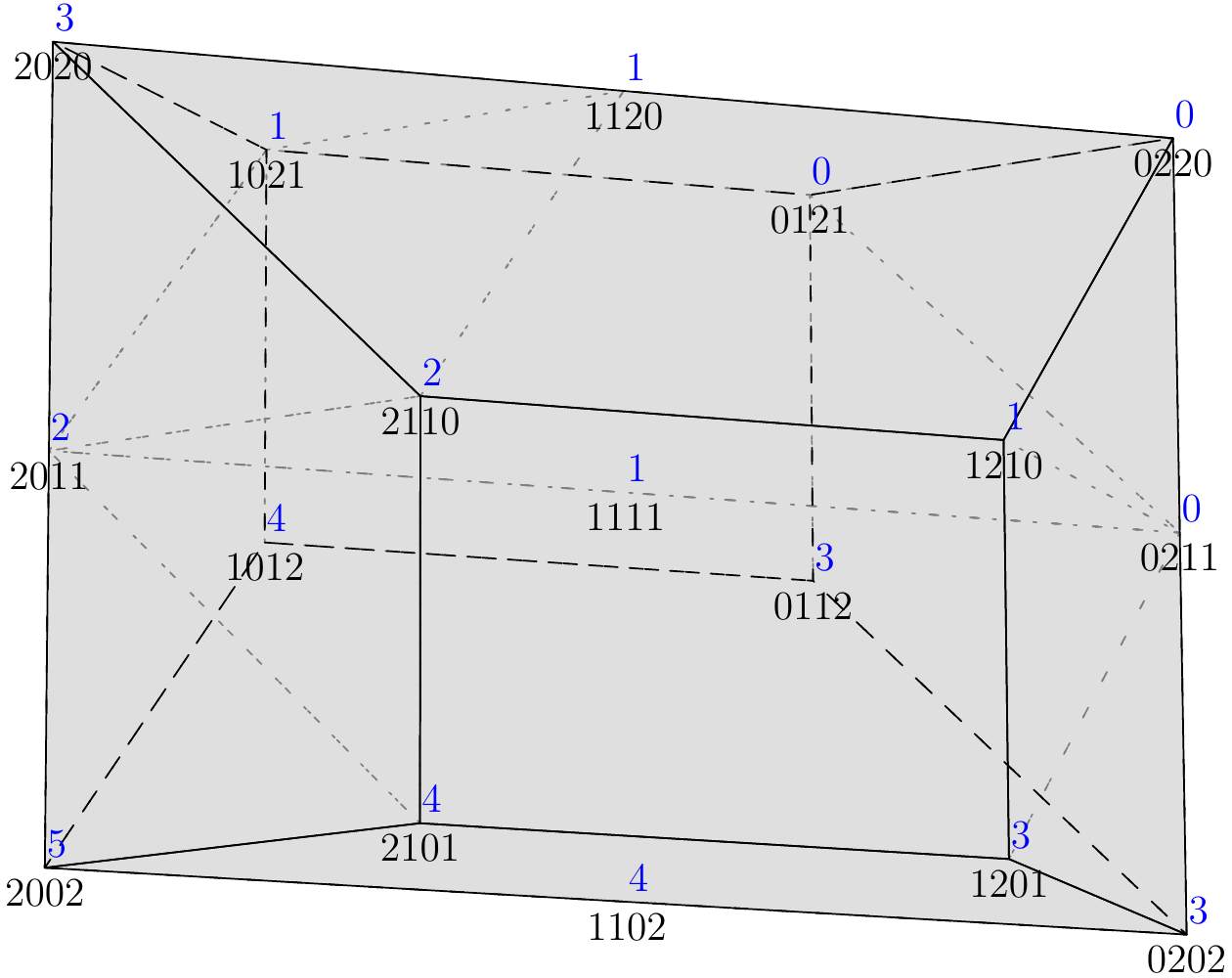}

\includegraphics[width=0.181818\textwidth]{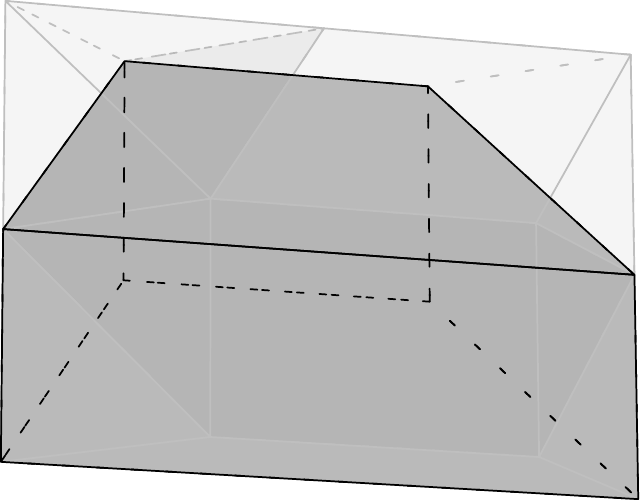} 
\includegraphics[width=0.181818\textwidth]{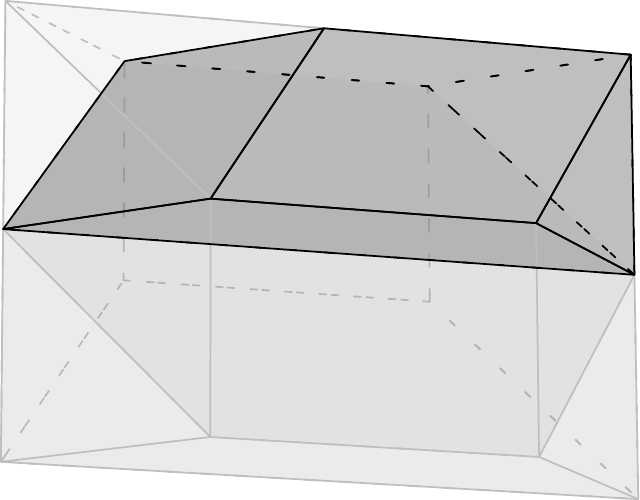}
\includegraphics[width=0.181818\textwidth]{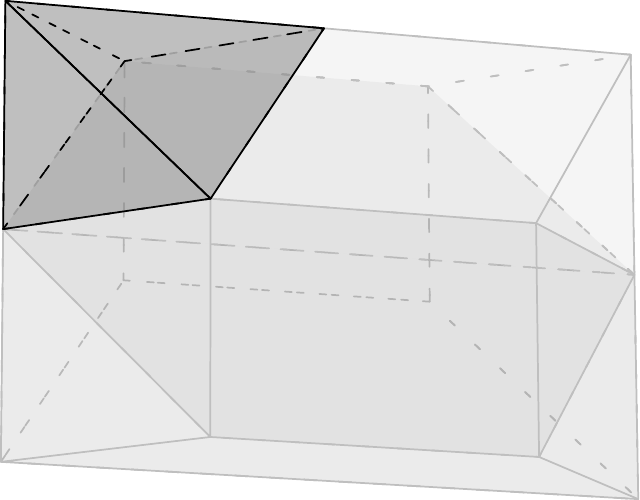}
\includegraphics[width=0.181818\textwidth]{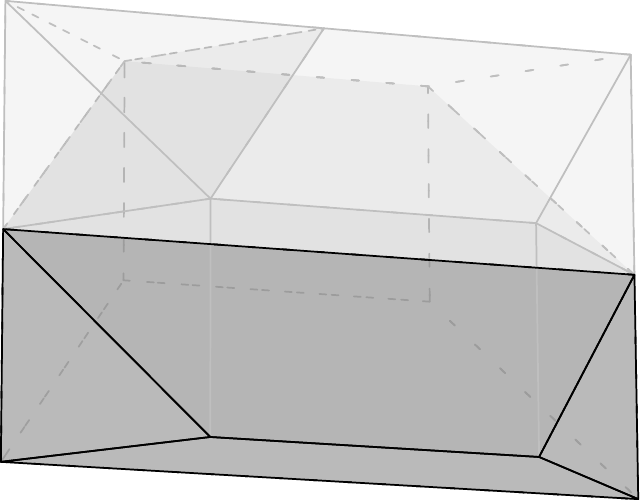}
\includegraphics[width=0.181818\textwidth]{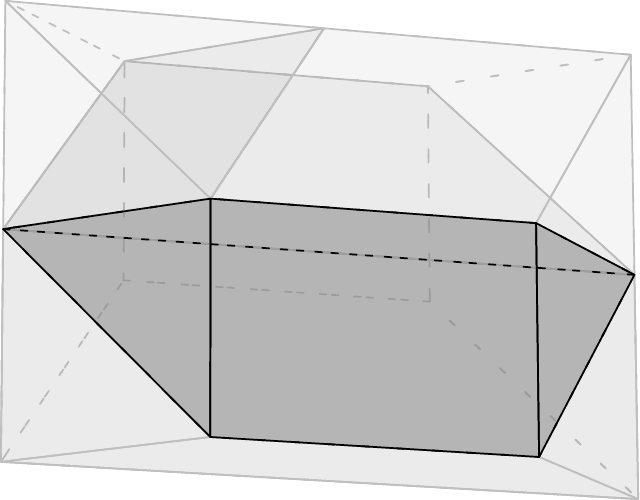}
}
\caption{The Chow subdivision of Example~\ref{ex:Chow subdivision}.  
Top: coordinates of points (black) and lifting heights (blue).
Bottom: the pieces.}
\label{fig:Chow subdivision}
\end{figure}

Perhaps the simplest varieties not among our opening examples
are conic curves in~$\bb P^3$.  
Let $\bbK=\bb C\lB t\rB$, and let
$\cl X\subseteq\bb P^3$ be the conic defined by the ideal 
$$(tx-y+z-t^3w, yz+tz^2+t^2yw-zw+(t^3-t^7)w^2)$$
where $(x:y:z:w)$ are coordinates on~$\bb P^3$.  
The Chow form of $\cl X$ can be computed by the algorithm of~\cite[Section~3.1]{Dalbec}.  It is 
{\footnotesize
\begin{multline*}
  (2t^7 + t^6 + t^5 - t^3)[zw][yw] + (t^7 + t^5 - t^3)[yw]^2 
+ (2t^4 + t^3 + t^2 - 1)[zw][yz] + (-t^3 + t^2 - 1)[yw][yz] \\
+ (-t - 1)[yz]^2 + (2t^8 - t^4)[zw][xw] + (2t^8 + t^6 - 2t^4)[yw][xw] 
+ (t^9 - t^5)[xw]^2 + (2t^5 - t)[zw][xz] \\+ (-t^4 + t^3 - 2t)[yw][xz] 
+ (-2t^2 - t)[yz][xz] + (-t^2)[xw][xz] + (-t^3)[xz]^2 
+ (-t^4 - 2t^3 + t)[zw][xy] \\+ (-t^3)[yw][xy] + t[yz][xy] 
+ (-t^4)[xw][xy] + t^2[xz][xy].
\end{multline*}
}
The Chow subdivision $\ChowSub(\cl X)$ is the regular subdivision induced by
the valuations of these coefficients.  It is a 3-polytope subdivided into 5 pieces,
depicted in Figure~\ref{fig:Chow subdivision}.  The polytope $\Chow(\cl X)$
of which it is a subdivision is an octahedron with two opposite corners truncated
(it is not the whole octahedron, which is the generic Chow polytope for conics in $\bb P^3$).
\end{example} 


Chow varieties and polytopes can also be defined for cycles on some more general spaces.
For this we of course suspend the assumptions of Case~\ref{case:projective}.
The groundwork for this construction is done in~\cite[Section~I.3]{Kollar},
and it's also treated in~\cite{KatzToolkit}.

\begin{case}\label{case:toric}
Let $\cl\iota:\cl Y\subseteq\bb P^{n-1}$ be a projective toric variety
with torus $T$, included $T$-equivariantly in $\bb P^{n-1}$.  All our Chow constructions
depend on $\cl\iota$, not merely $\cl Y$ alone.
Let $\Delta$ be the fan associated to~$Y$, and $N_\R$ its underlying vector space,
so that the fan structure defined on $N_\R$ by its 
intersections with cones of the fan $\mc L$ of $\bb P^{n-1}$ is equal to $\Delta$.
The inclusion $\cl\iota$ corresponds to a linear inclusion $\iota:N_\R\hookrightarrow\R/(1,\ldots,1)$, 
whose image we identify with~$N_\R$, turning tropical cycles in~$N_\R$ into tropical
cycles in~$\R/(1,\ldots,1)$.  These identifications are compatible with the corresponding classical ones.

Cycles in $\cl Y$ and $Y$ inherit a degree via $\cl\iota$ and $\iota$ respectively.  
For any given dimension $d-1$ and degree $r$, the Chow variety of dimension $d-1$ degree $r$
cycles for~$\cl\iota$ is defined as 
the subvariety of $\Gr(d,n,r)$ whose points represent cycles in $\cl Y$.  
By Theorem~\ref{r:slice-project}, the transpose $\iota^{\mathrm T}$ projects the simplex 
$\conv\{\eM{i}:i\in[n]\}$ 
onto a polytope $Q$ with $\mc N(Q)=\Delta$; this is the polytope associated to the ample divisor
$\cl\iota^*\mc O(1)$.  We define the Chow polytope and subdivision using the same projection.
For $\cl X\subseteq\cl Y$ a cycle, we define $\Chowi(\cl X) = \iota^{\mathrm T}\Chow(\cl X)$.
Similarly, if $\Pi\subseteq M_\R\times\R$ is the polytope determining the regular subdivision 
$\ChowSub(\cl X)$, then we define $\ChowiSub(\cl X)$ to be the regular subdivision of
$(\iota^{\mathrm T}\times\id_\R)(\Pi)$.
\end{case}

Returning to $\bb P^{n-1}$ as ambient variety,
Theorem~2.2 of~\cite{DFS} provides a procedure that determines 
the polytope~$\Chow(\cl X)$ given a fan tropical variety~$X=\Trop\cl X$.
That procedure is the constant-coefficient case of the next theorem, 
Theorem~\ref{r:DFS2.2}, which can be interpred as justifying our definition
of the Chow subdivision.
Theorem~\ref{r:DFS2.2} determines
$\ChowSub(\cl X)$ for $X=\Trop\cl X$ not necessarily a fan, by identifying the
regions of the complement of~$\mc N^1(\ChowSub(\cl X))$ and the 
vertex of~$\ChowSub(\cl X)$ each of these regions is dual to.

\begin{theorem}\label{r:DFS2.2}
Let $\dim X=d-1$.
Let $u\in N_\R$ be a linear functional such that 
$\face_u\ChowSub(\cl X)$ is a vertex of~$\ChowSub(\cl X)$.  
Then 
$$\inn_u R_{\cl X} = \prod_{J\in\binom{[n]}{n-d}}
[J]^{\deg([u+C_J]\cdot X)},$$
i.e.\
\begin{equation}\label{eq:DFS2.2 Chow}
\vertex_u\ChowSub(\cl X) = \sum_{J\in\binom{[n]}{n-d}}
\deg([u+C_J]\cdot X) \eM{J}.
\end{equation}
\end{theorem}

Recall that $C_J = \Rp\{\eN{j} : j\in J\}$.
The condition that $\face_u\ChowSub(\cl X)$ be a vertex is the
genericity condition necessary for the set-theoretic intersection 
$(u+C_J)\cap X$ to be a finite set of points.

The constant-coefficient case of Theorem~\ref{r:DFS2.2} is known as
{\em ray-shooting}, and the general case as {\em orthant-shooting},
since the positions of the vertices of $\ChowSub(\cl X)$ are read off
from intersection numbers of $X$ and orthants $C_J$ shot from the point~$u$.

\begin{example}
\begin{figure}
\centering{
\includegraphics{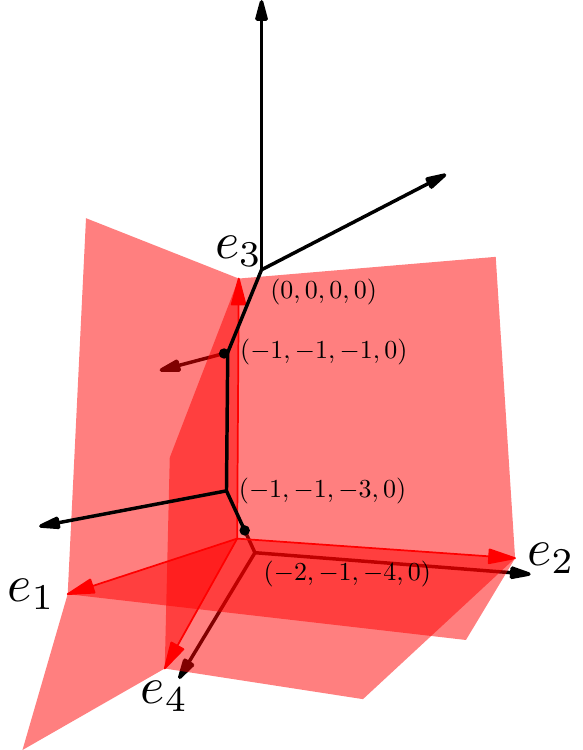}
}
\caption{Identifying a vertex of a Chow subdivision by
Theorem~\ref{r:DFS2.2}.  Coordinates of vertices of the curve
are given in $\R^4/(1,1,1,1)$.}
\label{fig:orthant-shooting}
\end{figure}

Let $\cl X$ be the conic curve of Example~\ref{ex:Chow subdivision}.
The black curve in Figure~\ref{fig:orthant-shooting} is $X=\Trop\cl X$.
Arbitrarily choosing the cone point of the red tropical plane to be $u\in N_\R$, 
we see that there are two intersection points among the
various $[u+C_J]\cdot X$, the two points marked as black dots.
Each has multiplicity 1, and they occur one each for $J=\{1,3\}$ and $J=\{2,3\}$.
Accordingly $\eM{\{1,3\}}+\eM{\{2,3\}}=(1,1,2,0)$ is the corresponding
vertex of~$\ChowSub(X)$ (compare Figure~\ref{fig:Chow subdivision}).
\end{example}

Theorem~\ref{r:DFS2.2} is proved in the literature, in a few pieces.
The second assertion, orthant-shooting in the narrow sense,
for arbitrary valued fields is Theorem~10.1 of~\cite{KatzToolkit}.
The first assertion, describing initial forms in the Chow form, is essentially
Theorem~2.6 of~\cite{KSZ}.  This is stated in the trivial valuation case
but of course extends to arbitrary valuations with our machinery of regular subdivisions
in one dimension higher.  The connection of that result with orthant shooting
is as outlined in Section~5.4 of~\cite{STY}.


For arbitrary ambient varieties, the second statement of Theorem~\ref{r:DFS2.2}
takes the following form.  (The analogue of the first statement is in \cite[Section~10]{KatzToolkit}.)
If $N_\R$ is a linear subspace of~$N'_\R$, and $X$ and~$Y$ are tropical cycles in~$N_\R$,
by $X\mathbin{\cdot_{N_\R}}Y$ we mean the stable tropical intersection taken in $N_\R$,
i.e.\ where the displacement in the fan displacement rule is restricted to $N_\R$.

\begin{corollary}
With setup as in Case~\ref{case:toric}, let $\dim\cl X=d-1$ let $u\in N_\R$ be such that 
$\face_u\ChowiSub(\cl X)$ is a vertex of $\ChowiSub(\cl X)$.  Then
\begin{equation}\label{eq:toric Chow}
\vertex_u\ChowiSub(\cl X) = \sum_{F}
\deg\big([u+\normal_{N_\R}(F)]\mathbin{\cdot_{N_\R}} X\big)\, m(F)
\end{equation}
where $F$ runs over faces of~$Q$ of dimension $\dim X$, and
$m(F) = d\int_F x\,\mathrm dx\in \R^n$.
\end{corollary}

\begin{proof}
By definition $\ChowiSub(\cl X)$ is the image of $\ChowSub(\cl X)$ under $\iota^{\mathrm T}$.
For any tropical cycles $X\subseteq Y$ and $Z$, we have that 
$(Z\cdot Y)\mathbin{\cdot_Y}X = Z\cdot X$ (in treatments
such as~\cite{AR}, which develop tropical cycles as zero loci of collections
of rational functions and intersection as restriction of rational functions,
this is immediate).  Using this in~\eqref{eq:DFS2.2 Chow} gives
\begin{align*}
\vertex_u\ChowiSub(\cl X) &= \sum_{J\in\binom{[n]}{n-d}}
\deg\big(([u+C_J]\cdot [N_\R])\mathbin{\cdot_{N_\R}} X\big)\, \iota^{\mathrm T}(\eM{J})
\\&=\sum_{J\in\binom{[n]}{n-d}}
\deg\big(([u]\msc([C_J]\cdot [N_\R]))\mathbin{\cdot_{N_\R}} X\big)\, \iota^{\mathrm T}(\eM{J}).
\end{align*}
By Theorem~\ref{r:slice-project},
the sum of these $[C_J]\cdot[N_\R]$ is $\mc N^{d-1} Q$, with the natural tropical weights.
A cone of $\mc N^{d-1} Q$ may arise from multiple cones $C_J$.  
For each dimension $d-1$ face $F$ of~$Q$,
consider the images of those vertices of $\conv\{\eM{i}:i\in[n]\}$ mapped onto it by $\iota^{\mathrm T}$.
Take a triangulation of $F$ using these images, {\em coherent} in the sense of~\cite[Section~4]{McMullen2},
and suppose the simplices used are the $S_J := \conv\{\iota^{\mathrm T}\eM{i}:i\in J\}$ 
for $J\in T(F)\subseteq\binom{[n]}{d}$.  Then we have 
$$
\vertex_u\ChowiSub(\cl X) = \sum_{F} \deg([u + \normal_{N_\R}(F)]\mathbin{\cdot_{N_\R}} X)\, \vol(F)
\sum_{J\in T(F)}\iota^{\mathrm T}(\eM{J}).
$$
But for each $J$ we have $\int_{S_J} x\,\mathrm dx = \vol(S_J)\,\iota^{\mathrm T}(e^J)/(d)$.
Summing this integral over all the simplices in~$F$ yields~\eqref{eq:toric Chow}.
\end{proof}

\section{From tropical variety to Chow polytope}\label{sec:trChow}

Henceforth $d\leq n$ will be a fixed integer, and
$\cl X$ will be a $(d-1)$\/-dimensional subvariety of the ambient toric variety,
which is mostly $\bb P^{n-1}$.  

As explained in~\cite{KSZ},
the torus $(\bbK^\ast)^n$ acts on the Hilbert scheme $\Hilb(\bb P^{n-1})$
in the fashion induced from its action on~$\bb P^{n-1}$, and the map 
$\Hilb(\bb P^{n-1})\to\Gr(d,n,r)$ sending each ideal to the
corresponding cycle is $(\bbK^\ast)^n$-equivariant.
This implies that deformations in $\Hilb(\bb P^{n-1})$
determine those in $\Gr(d,n,r)$: if $u,u'\in N_\R$ are such that 
$\inn_u\mc I(\cl X) = \inn_{u'}\mc I(\cl X)$, where $\mc I$ denotes the defining ideal, 
then also $\inn_u R_{\cl X}=\inn_{u'} R_{\cl X}$.  Accordingly each
initial ideal  of~$\mc I(\cl X)$ determines a face of $\Chow(\cl X)$, so that 
the Gr\"obner fan of $\cl X$ is a refinement of the normal fan 
of~$\Chow(\cl X)$. 

The standard construction of the tropical variety $X$ via initial ideals \cite[Theorem~2.6]{RST}
shows that $X$ is a subfan of the Gr\"obner fan.
But in fact $X$ is a subfan of the coarser fan $\mc N(\Chow(\cl X))$,
since the normal cone of a face $\face_u\Chow(\cl X)$ appears 
in~$X$ if and only if $\cl X$ meets the maximal torus 
$(\bbK^\ast)^n/\bbK^\ast\subseteq\bb P^{n-1}$, and whether this happens is
determined by the cycle associated to~$\cl X$.  
The analogue of this holds in the non-fan case as well.
This reflects the principle that the information encoded in the Hilbert scheme but not 
in the Chow variety pertains essentially to nonreduced structure, while
tropical varieties have no notion of embedded components and only
multiplicities standing in for full-dimensional non-reduced structure.

The machinery of Section~\ref{sec:msc} allows us to give a lean combinatorial characterisation 
of the Chow subdivision in terms of Theorem~\ref{r:DFS2.2}.

\begin{maintheorem}\label{r:1}
For projective tropical varieties, we have
\begin{equation*}
\mc N^1(\ChowSub(\cl X)) = X \msc {\mc L_{n-d-1}}^\refl.
\end{equation*}
In general, with the notation of Case~\ref{case:toric},
$$
\mc N^1(\ChowiSub(\cl X)) = X \msc {\mc N^d(Q)}^\refl.
$$
\end{maintheorem}
\noindent To reiterate:
{\em Let $\cl X$ be a $(d-1)$-cycle in a projective tropical variety $\cl Y$, 
and let $X=\Trop\cl X$.
Then the codimension 1 part of the normal subdivision to the Chow subdivision of~$\cl X$
is the stable Minkowski sum of~$X$ and the reflection of the codimension $d$ skeleton 
of the fan of $\cl Y$} (with its natural weights under the embedding).  
In the projective case, the second summand is the reflected linear space ${\mc L_{n-d-1}}^\refl$.
By Theorem~\ref{r:weights on polytopes}(b), this uniquely determines $\ChowSub(\cl X)$
in terms of~$X$, up to translation and adding a constant to the vertex heights.

Theorem~\ref{r:1} should be taken as providing the extension
of the notion of Chow polytope (via its normal fan) to tropical varieties.
\begin{definition}
Let the {\em Chow map} $ch$ for projective space be the map taking a tropical cycle $X$ of dimension $d$
to its ({\em tropical}) {\em Chow hypersurface}, the cycle $ch(X) = X \msc {\mc L_{n-d-1}}^\refl$.

More generally, for an ambient projective toric variety $\cl\iota:\cl Y\to\bb P^{n-1}$, let the
Chow map $ch_{\cl\iota}$ be given by $ch_{\cl\iota}(X) = X\msc\mc N^d(Q)^\refl.$
\end{definition}
\noindent The dimension of $ch(X)$ is $(d-1)+(n-d-1)=n-2$, so its codimension is~1.
Indeed $ch$ is a linear map $ Z_{d-1}\to Z^1$.  Likewise $ch_\iota$ is a linear map 
$Z_{d-1}\to Z_{\dim Y-1}$.

\begin{remark}\label{rem:cl vs trop}
In the projective case, the support of $ch(X)$ is precisely the set of points $u\in N_\R$ such that
a tropical $(n-d-1)$-plane centered at $u$ meets $X$.  This is very reminiscent of the 
classical construction of the Chow form in Remark~\ref{rem:meet}, which uses classical
$(n-d-1)$-planes meeting~$\cl X$.  The most significant difference between the two
constructions is that the classical {\em Chow hypersurface}
lies in $\Gr(n-d,n)$, where it is the zero locus of the Chow form~$R_{\cl X}$.  
By contrast our tropical Chow hypersurface $ch(X)$ 
lies in the tropical torus $(\bbK^\ast)^n/\bbK$, in the same space as~$X$.  
One might think of this as reflecting the presence in tropical projective geometry
of a single canonical nondegenerate linear space ${\mc L_e}$ of each dimension, something with
no classical analogue.

Following the classical construction more closely, one could associate to~$\cl X$
a hypersurface $Y$ in~$\Trop\Gr(n-d,n)$, namely
the tropicalisation of the ideal generated by $R_{\cl X}$ and the Pl\"ucker relations.  
The torus action $(\bbK^\ast)^n/\bbK^\ast\curvearrowright\Gr(n-d,n)$
tropicalises to an action of $N_\R$ on~$\Trop\Gr(n-d,n)$ by translation,
i.e.\ an $(n-1)$-dimensional lineality space.  
Denote by $N_\R+0$ the orbit of the origin in~$\Trop\Gr(n-d,n)$;
this is the parameter space for tropical linear spaces in $N_\R$ that 
are translates of~$\mc L_{n-d-1}$.  Then we have $ch(X) = Y\cap (N_\R+0)$.  
\end{remark}

Lemma~\ref{r:ch deg} is also seen to be about Chow hypersurfaces, in which context it says
$$\deg ch(X) = \codim X\deg X.$$
This should be compared to the fact that
the Chow form of a cycle $\cl X$ in $\Gr(d,n,r)$ is of degree $r=\deg\cl X$ in~$\bb K[\Gr(n-d,n)]$,
and this ring is generated by brackets in $n-d=\codim\cl X$ letters.

\begin{proof}[Proof of Theorem~\ref{r:1}]
We begin in the projective Case~\ref{case:projective}.
Given a regular subdivision $T$ of lattice polytopes in~$M$ induced by~$\Pi$,
its support function $V_T:u\mapsto\face_u T$ is a piecewise linear function
whose domains of linearity are~$\mc N^0(T,\Pi)$.
We can view $V_T$ as an element of~$(Z^\unbal)^0\otimes M$.

We take a linear map $\delta:( Z^\unbal)^0\otimes M\to ( Z^\unbal)^1$ 
such that $\delta(V_T) = \mc N^1(T,\Pi)\in Z^1$ for any regular subdivision $T$.
The restriction of $\delta$ to the linear span of all support functions is
a canonical map $\delta'$, which has been constructed as the map 
from Cartier divisors supported on $\mc N(T,\Pi)$ to Weil divisors on $\mc N(T,\Pi)$ in
the framework of~\cite{AR}, or as the map from piecewise polynomials to Minkowski weights 
given by equivariant localisation in~\cite{KP}.
Roughly, $\delta'(V)$ is the codimension 1 tropical cycle whose multiplicity at
a facet $\tau$ records the difference of the values taken by~$V$ on either side of~$\tau$.
We can take $\delta$ as any linear map extending $\delta'$ such that $\delta(V)$
still only depends on $V$ locally;
our only purpose in making this extension is to allow formal manipulations using unbalanced cycles.

%
%

Let $V=V_{\ChowSub(\cl X)}$, and
write $X=\sum_{\sigma\in\Sigma} m_\sigma[\sigma]$.
Expanding \eqref{eq:DFS2.2 Chow}
in terms of this sum, the value of~$V$ at $u\in N_\R$ is
$$\sum_{\sigma\in\Sigma} m_\sigma
\sum_{J\in\binom{[n]}{n-d}} \deg([\sigma]\cdot[u+C_J])\eM{J}.$$
The intersection $[\sigma]\cdot[u+C_J]$ is zero if $u\not\in\sigma-C_J$,
and if $u\in\sigma-C_J$
it is one point with multiplicity $\mu^\bullet_{\sigma,C_J}$.
So
$$V = \sum_{\sigma\in\Sigma} m_\sigma
\sum_{J\in\binom{[n]}{n-d}} \mu^\bullet_{\sigma,C_J}[\sigma-C_J]\otimes \eM{J}.$$
Let $V_\sigma$ be the inner sum here, so that $V=\sum_{\sigma\in\Sigma} m_\sigma V_\sigma$.
Then 
$$\delta(V_\sigma) = \sum_{J\in\binom{[n]}{n-d}} \mu^\bullet_{\sigma,C_J}
\sum_{\stackbounds{$\tau$ a facet}{of $\sigma-C_J$}} \delta([\tau]\otimes \eM{J}).$$
Here, if $\tau$ is a facet of form $\sigma'-C_J$ for~$\sigma'$ a
facet of~$\sigma$, then $\eM{J}\in\R\tau$ so $\delta([\tau]\otimes \eM{J})=0$ and the
$\tau$ term vanishes.  Otherwise $\tau$ has the form 
$\sigma-C_{J'}$ where $J'=J\setminus\{j\}$ for some $j\in J$.  Regrouping
the sum by~$J'$ gives
\begin{equation}\label{eq:pd0fs}
\delta(V_\sigma) = \sum_{J'\in\binom{[n]}{n-d-1}} 
\left(\sum_{j\in[n]\setminus J'}
\mu^\bullet_{\sigma,C_{J'\cup\{j\}}}\delta([\sigma-C_{J'}]\otimes \eM{j})\right)
\end{equation}
where again we have omitted the terms $\delta([\sigma-C_{J'}]\otimes \eM{J'})=0$.  Now,
if $j\not\in J'$ then 
\begin{align*}
\mu^\bullet_{\sigma,C_{J'\cup\{j\}}}
&= \mu_{\sigma,C_{J'\cup\{j\}}}
\\&= [N_{\sigma+C_{J'\cup\{j\}}} : N_\sigma + N_{C_{J'\cup\{j\}}}]
\\&= [N : N_\sigma + N_{C_{J'}} + \bb Z\eN{j}]
\\&= [N : N_{\sigma+C_{J'}} + \bb Z\eN{j}]
        [N_{\sigma+C_{J'}} + \bb Z\eN{j} : N_\sigma + N_{C_{J'}} + \bb Z\eN{i}]
\\&= [N : N_{\sigma+C_{J'}} + \bb Z\eN{j}][N_{\sigma+C_{J'}} : N_\sigma + N_{C_{J'}}]
\\&= \langle \eN{j},p\rangle \mu^\msc_{\sigma,C_{J'}}
\end{align*}
where $p$ is the first nonzero lattice point in the appropriate direction 
on a line in~$M_\R$ normal to~$\sigma+C_{J'}$.  Then the components of~$p$
are the minors of a matrix of lattice generators for~$\sigma+C_{J'}$ 
by Cramer's rule, and the last equality is a row expansion of 
the determinant computing $\mu^\msc_{\sigma,C_{J'}}$.  
If $j\in J'$ then 
$\mu^\bullet_{\sigma,C_{J'\cup\{j\}}}=0=\langle \eN{j},p\rangle \mu^\msc_{\sigma,C_{J'}}$
also.  So it's innocuous to let the inner sum in~\eqref{eq:pd0fs}
run over all~$j\in[n]$, and we get
\begin{align*}
\delta(V_\sigma) &= \sum_{J'\in\binom{[n]}{n-d-1}} 
\left(\sum_{j\in[n]}
\mu^\msc_{\sigma,C_{J'}}\langle \eN{j},p\rangle\delta([\sigma-C_{J'}]\otimes \eM{j})\right)
\\&= \sum_{J'\in\binom{[n]}{n-d-1}}\mu^\msc_{\sigma,C_{J'}}
\delta([\sigma-C_{J'}]\otimes p)
\\&= \sum_{J'\in\binom{[n]}{n-d-1}}\mu^\msc_{\sigma,C_{J'}}[\sigma-C_{J'}]
\\&= ([\sigma] \msc {\mc L_{n-d-1}}^\refl).
\end{align*}
We conclude that 
$$\mc N^1(\ChowSub(\cl X))=\delta(V)=\sum_\sigma m_\sigma ([\sigma] \msc {\mc L_{n-d-1}}^\refl)
= X\msc {\mc L_{n-d-1}}^\refl.$$  

Finally we handle the case of arbitrary ambient variety.  We have that
${\mc L_{n-d-1}}$ is the codimension~$d$ skeleton of the simplex $S:=\conv\{e^i:i\in[n]\}$.
Then
\begin{align*}
\mc N^1(\ChowiSub(\cl X)) 
&= \mc N^1(\iota^{\mathrm T}\ChowSub(\cl X))
\\&= Y\cdot \mc N^1(\ChowSub(X)) &\mbox{by Theorem~\ref{r:slice-project}}
\\&= Y\cdot (X\msc\mc N^d(S))
\\&= X\msc (Y\cdot\mc N^d(S)) &\mbox{by Lemma~\ref{r:slicing commutes}}
\\&= X\msc\mc N^d(Q) &\mbox{by Theorem~\ref{r:slice-project}.\qedhere}
\end{align*}
\end{proof}

\section{Linear spaces}\label{sec:linear}

A {\em matroid subdivision} (of rank $r$) is a regular subdivision of a matroid
polytope (of rank $r$) all of whose facets are matroid polytopes,
i.e.\ polytopes of the form $\Poly(M)$ defined in~\eqref{eq:matroid polytope}.
The {\em hypersimplex} $\Delta(r,n)$ is the polytope
$\conv\{\eM{J} : J\in\binom{[n]}r\}$.  The vertices of
a rank~$r$ matroid polytope are a subset of those of~$\Delta(r,n)$.
We have the following polytopal characterisation of matroid polytopes due to
Gelfand, Goresky, MacPherson, and Serganova.

\begin{theorem}[\cite{GGMS}]\label{r:GGMS} A polytope $\Pi\subseteq \R^n$
is a matroid polytope if and only if $\Pi\subseteq[0,1]^n$ 
and each edge of~$\Pi$ is a parallel translate of $\eM{i}-\eM{j}$ for some $i,j$.
\end{theorem}

\begin{definition}\label{def:Bergman}
Given a regular matroid subdivision $\Sigma$, its {\em Bergman complex} $\Berg(\Sigma)$ 
and {\em co-Bergman complex} $\coBerg(\Sigma)$ are subcomplexes of $\mc N(\Sigma)$.
The face of $\mc N(\Sigma)$ normal to $F\in\Sigma$
\begin{itemize}
\item is a face of $\Berg(\Sigma)$ if and only if $F$ is the polytope of a loop-free matroid;
\item is a face of $\coBerg(\Sigma)$ if and only if $F$ is the polytope of a coloop-free matroid.
\end{itemize}
We make $\Berg(\Sigma)$ and $\coBerg(\Sigma)$ into tropical varieties 
by giving each facet multiplicity~1.
\end{definition}

The {\em Bergman fan}, the fan case of the Bergman complex, was introduced in~\cite{AK}
(where an object named the ``Bergman complex'' different to ours also appears).
Bergman complexes are much used in tropical geometry, on account of the 
following standard definition, appearing for instance in~\cite{Speyer}.
\begin{definition}\label{def:TLS alt}
A {\em tropical linear space} is the Bergman complex of a regular matroid subdivision.  
\end{definition}
In the context of Chow polytopes it is the co-Bergman complex 
rather than the Bergman complex that arises naturally,  
on account of the duality mentioned in Example~\ref{ex:motivating2}(2).
Observe that the co-Bergman complex of a matroid subdivision
is a reflection of the Bergman complex of the dual matroid subdivision;
in particular any Bergman complex is a co-Bergman complex and vice versa.

Since there is a good notion of tropical degree (Definition~\ref{def:degree}), 
the following alternative definition seems natural.
\begin{definition}\label{def:TLS}
A {\em tropical linear space} is a tropical variety of degree~1.
\end{definition}

\begin{theorem}\label{r:TLS}
Definitions \ref{def:TLS alt} and~\ref{def:TLS} are equivalent.
\end{theorem}

The equivalence in Theorem~\ref{r:TLS} was noted by Mikhalkin, Sturmfels, and Ziegler
and recorded in~\cite{IMS}, but no proof was provided.   
One implication, that Bergman complexes of matroids have degree~1,
follows from Proposition~3.1 of~\cite{Speyer}, which implies that 
the tropical stable intersection of a $(d-1)$-dimensional Bergman complex of a matroid subdivision 
with $\mc L_{n-d}$ (the Bergman complex of a uniform matroid) is a 
0-dimensional Bergman complex, i.e.\ a point with multiplicity 1.
Thus it remains to prove that degree~1 tropical varieties are (co-)Bergman complexes.  
In fact, let $X\subseteq N_\R$ be a degree~1 tropical variety of dimension~$d-1$. 
We will show
\begin{enumerate}
\item The regular subdivision $\Sigma$ such that 
$ch(X)=\mc N^1(\Sigma)$ is dual to a matroid subdivision of rank $n-d$. 
\smallskip
\item We have $X=\coBerg(\Sigma)$.
\end{enumerate}

Tropical varieties have an analogue of B\'ezout's theorem.
See for instance Theorem~9.16 of~\cite{AR}, which however only proves
equality under genericity assumptions, not the inequality below.
We will only need the theorem in the case that the varieties being intersected have degree~1.

\begin{theorem}[Tropical B\'ezout's theorem]\label{r:bezout}
Let $X$ and~$Y$ be tropical varieties of complementary dimensions.  We have
$\deg(X\cdot Y)\leq \deg X\deg Y$, and equality is attained
if $X$ and~$Y$ are of sufficiently generic combinatorial type.  
\end{theorem}


\begin{lemma}\label{r:TLSs have no neg rays}
If a tropical variety $X$ of degree~1 contains a ray in direction $-\eN{i}$ for $i\in[n]$,
then $-\eN{i}$ is contained in the lineality space of~$X$. 
\end{lemma}

\begin{proof}
Consider the set 
$$Y=\{u\in N_\R: \mbox{$u-a\eN{i}\in X$ for $a\gg0$}\}.$$
By assumption on~$X$, $Y$ is nonempty.
This $Y$ is the underlying set of a polyhedral complex; make it into a cycle by giving
each facet multiplicity~1.  In fact, $Y$ is a tropical variety, 
as any face $\tau$ of~$Y$ corresponds to a face $\sigma$ of~$X$ 
such that $\tau=\sigma+\R \eN{i}$, and so $Y$ inherits balancing from~$X$.  
Also $\dim Y=\dim X=:d-1$.  
Since $Y$ is effective, some translate and therefore any translate of~$\mc L_{n-d-1}$ 
intersects $Y$ stably in at least one point.

Suppose $X$ had a facet $\sigma$ whose linear span didn't contain $-\eN{i}$. 
Then there is some translate $[u]\msc\mc L_{n-d-1}$ which intersects $\relint\sigma$,
with the intersection lying on a face $u+C_J$ of~$[u]\msc\mc L_{n-d-1}$ with $i\in J$.
Given this translate, any other translate $[u-a\eN{i}]\msc\mc L_{n-d-1}$ with $a\geq0$ 
will intersect~$X$ transversely in the same point of~$\relint\sigma$.  
For $a$ sufficiently large, one of the 
points of $Y\cdot([u-a\eN{i}]\msc\mc L_{n-d-1})$ lies in~$X$, providing 
a second intersection point of $X$ and~$[u-a\eN{i}]\msc\mc L_{n-d-1}$.
By B\'ezout's theorem this contradicts the assumption that $\deg X=1$.
\end{proof}

\begin{proof}[Proof of Theorem~\ref{r:TLS}]
{\em To} (1).
Suppose $l\subseteq N_\R$ is a classical line in any direction~$\eN{J}$, $J\subseteq[n]$.
By Lemma~\ref{r:two pairings} and Theorem~\ref{r:bezout} we have
\begin{equation}\label{eq:chX cdot l}
\deg(ch(X)\cdot [l]) = \deg((X\msc {\mc L_{(n-d-1)}}^\refl)\cdot [l])
= \deg((\mc L_{(n-d-1)}\msc [l])\cdot X)\leq 1
\end{equation}
because $\mc L_{(n-d-1)}\msc [l]$ is a degree~1 tropical variety.
Since intersection multiplicities are positive,
if $l$ intersects a facet $\sigma$ of~$ch(X)$
then the multiplicity of the intersection is $\mu^\bullet_{\sigma,l}=1$.

Let $\sigma$ be a facet of~$ch(X)$, and $l$ a line in
direction~$\eN{J}$ intersecting it.
Then $\mu^\bullet_{\sigma,l} = \langle m, \eN{J}\rangle$ where $m\in M_\R$ is 
the difference of the endpoints of the edge of~$\Sigma$ dual to $\sigma$.
Then $m$ is the product of a primitive normal vector to~$\sigma$ and the multiplicity $m_\sigma$.  
The positive components of~$m$ cannot have sum $k\geq2$, or else,
for a suitable choice of~$J$, 
we would achieve $\mu^\bullet_{\sigma,l}=\langle m,\eN{J}\rangle=k$.  
Since $m$ is nonzero and normal to~$(1,\ldots,1)$
we must have $m=\eM{i}-\eM{j}$ for some $i\neq j\in[n]$.  
It follows that each edge of~$\Sigma$ is a parallel translate of some~$\eM{i}-\eM{j}$.  

Furthermore, let $l\subseteq N_\R$ be a line in direction $\eN{i}$, for $i\in[n]$.
The vertices of~$\Sigma$ attained as $\face_u\Sigma$ for some $u\in l$
are in bijection with the connected components of the complement of~$ch(X)$.
So there are at most two of these vertices, 
and if there are two, say $m_0$ and~$m_1$, we have $\langle m_1-m_0, \eN{i}\rangle=1$.  But among
the vertices $\face_u\Sigma$ for $u\in l$ are vertices $m$ minimising and maximising
the pairing $\langle m,\eN{i}\rangle$.  Therefore, the projection of $\Sigma$
to the $i$\/th coordinate axis has length either 0 or~1.

For the remainder of the proof we fix a particular translation representative
of~$\Sigma$, namely the one whose projection onto the $i$\/th coordinate axis
is either the point $\{0\}$ or the interval $[0,1]$ for each $i\in[n].$ 
For this particular $\Sigma$, 
Theorem~\ref{r:GGMS} implies that $\Sigma$ is a matroid subdivision.

Let $r$ be the rank of the matroid subdivision $\Sigma$. 
Let $\eM{J}$ be one vertex of~$\Sigma$, so that $|J|\in\binom{[n]}r$, 
and let $u$ be a linear form with $\face_u\Sigma=\eM{J}$.  
Then, for any $i\in[n]\setminus J$ and any $a>0$, we have 
$\face_{u+a\eN{i}}\Sigma = \eM{J}$, since $\eM{J}\in\face_{\eN{i}}\Sigma$.  
On the other hand, for any $i\in J$ and sufficiently large $a\gg 0$,
we have $\face_{u+a\eN{i}}\Sigma\not\ni \eM{J}$, and indeed 
$\face_{u+a\eN{i}}\Sigma$ will contain some vertex $\eM{J'}$ with $i\not\in J'$,
whose existence is assured by our choice of translation representative for~$\Sigma$.
It follows that a ray $[u]\msc[R_{\geq0}\{\eN{i}\}]$ of $[u]\msc\mc L_1$
intersects $ch(X)$ if and only if $i\in J$.  Each intersection
must have multiplicity 1, so 
$$\deg(ch(X))=\deg(ch(X)\cap([u]\msc\mc L_1))=|J|=r.$$ 
But by Proposition~\ref{r:ch deg} we have that $\deg(ch(X)) = n-d$,
so $r=n-d$ as claimed.

\noindent {\em To} (2).
Fix some polyhedral complex structure on~$X$. 
Given any $u\in N_\R$ in the support of~$ch(X)$, 
its multiplicity is $ch(X)(u) = 1$, and therefore by positivity there is a unique choice of
a facet $\tau$ of~$X$ and $J\in\binom{[n]}{n-d-1}$ such that $u\in X-C_J$.
Write $J = J(u)$.  
On the other hand, $\Sigma$ has a canonical coarsest possible 
polyhedral complex structure, on account of being a normal complex.  
We claim that $J(u)$ is constant for
$u$ in the relative interior of each facet $\sigma$ of~$\Sigma$,
and thus we can write $J(\sigma) := J(u)$.
Suppose not. 
Consider the common boundary $\rho$ of two adjacent regions $\sigma_1,\sigma_2$
of~$\sigma$ on which $J(u)$ is constant.  Suppose $\sigma_1\subseteq \tau-C_{J_1}$.  
We have $\rho\subseteq\tau-C_K$ for $K\in\binom{[n]}{n-d-2}$.
There is a facet of~$\Sigma$ of form $\sigma_j\subseteq\tau-C_{K\cup k}$ 
incident to~$\rho$ for each $k\in[n]\setminus K$ such that $\eN{k}$ is not contained in
the affine hull of~$\tau$.  Since $\dim\tau=d-1$, and any $d$ of the $\eN{k}$ are independent
in~$N_\R$, there exist at most $d-1$ indices $k\in[n]$ such that $\eN{k}$ is not 
contained in the affine hull of~$\tau$, and hence at least
$$|[n]\setminus K| - (d-1) = 3$$
indices $k\in[n]$ yielding facets of $\Sigma$.  In particular $\sigma_1$ and~$\sigma_2$
cannot be the only $(d-1)$-dimensional regions
in~$\Sigma$ incident to~$\rho$, and this implies $\sigma$ cannot be a facet of~$\Sigma$,
contradiction.

Now, every facet $\sigma$ of~$ch(X)$ is normal to an edge of~$\Sigma$, say
$E_\sigma=\conv\{\eM{K}+\eM{j},\eM{K}+\eM{k}\}$ for $K\in\binom{[n]}{n-d-1}$.  Since $\Sigma\subseteq\Delta(n-d,n)$,
$\sigma$ must contain a translate of
the normal cone to~$E_\sigma$ in~$\mc N^1(\Delta(n-d,n))$, namely
$$\normal(E_\sigma) = \{u\in N_\R : u_j=u_k, \mbox{$u_i\leq u_j$ for $i\in K$}, 
\mbox{$u_i\geq u_j$ for $i\not\in K\cup\{j,k\}$}\}.$$ 
In particular $\sigma$ contains exactly
    $n-d-1$ rays in directions $-\eN{i}$, those with $i\in K$.

Let $R$ be the set of directions $-\eN{1},\ldots,-\eN{n}$.
Suppose for the moment that $X$ contains no lineality space in any direction $-\eN{i}$.  
We have that $\sigma\subseteq X\msc[-C_{J(\sigma)}]$.
By Lemma~\ref{r:TLSs have no neg rays},
$X$ contains no rays in directions in~$R$, so we must have 
that $J(\sigma)=K$ and $-C_{J(\sigma)}$ contains a ray in direction $-\eN{i}$ for all $i\in K$.  
Now consider any face $\rho$ of~$\sigma$ containing no rays in directions in~$R$.
Then we claim $\rho\in X$.  
If this weren't so, then there would be another face $\sigma'$ parallel to $\sigma$
and with $J(\sigma)=J(\sigma')$.  But the edge $E_\sigma$ is determined by $J(\sigma)=K$
and the normal direction to~$\sigma$, so $E_\sigma=E_{\sigma'}$, implying $\sigma=\sigma'$.
On the other hand, the relative interior of any face of~$\sigma$ 
containing a ray in direction~$R$ is disjoint from~$X$, since if $u$ is a point in such a face
there exists $v\in -C_{J(\sigma)}\setminus\{0\}$ such that $u-v\in X$.
So $X$ consists exactly of the faces of~$ch(X)$ containing no ray in a direction in~$R$.

If $X$ has a lineality space containing those $-\eN{j}$ with $j\in J$, then
let $X'$ be the pullback of~$X$ along a linear projection with kernel $\spann\{-\eN{j} : j\in J\}$.
Then we can repeat the last argument using $X'$, and we get that
$X$ consists exactly of the faces of~$ch(X)$ containing no ray in a direction in
$R\setminus\{-\eN{j}:j\in J\}$.

Now, a face $\normal(F)$ of~$\mc N(\Sigma)$ contains a ray in direction $-\eN{i}$
if and only if the linear functional $\langle m,-\eN{i}\rangle$ is constant on~$m\in F$ 
and equal to its maximum for~$m\in\Sigma$.  
The projection of $F$ to the $i$\/th coordinate axis is either $\{0\}$, $\{1\}$, or~$[0,1]$,
so $\normal(F)$ contains a ray in direction $-\eN{i}$ if and only if
the projection of~$F$ is~$\{1\}$, or the projection of~$F$ and of~$\Sigma$ are both~$\{0\}$.
Projections taking $\Sigma$ to~$\{0\}$ correspond to lineality directions in~$X$,
so we have that $X$ consists exactly of the faces of~$ch(X)$ which don't project
to~$\{1\}$ along any coordinate axis.  These are exactly the coloop-free faces.
\end{proof}



\section{The kernel of the Chow map}\label{sec:kernel}

In this section we will show that 
the Chow map $ch: Z_{d-1}\to Z^1$ has a nontrivial kernel.
This implies that there exist distinct tropical varieties with the same Chow polytope:
$Y$ and $X+Y$ will be a pair of such varieties for any nonzero $X\in\ker ch$, 
choosing $Y$ to be any effective 
tropical cycle such that $X+Y$ is also effective (for instance, let $Y$
be a sum of classical linear spaces containing the facets of~$X$
that have negative multiplicity).
Thus Chow subdivisions do not lie in a combinatorial bijection with general tropical
varieties, as was the case for our opening examples.  


There are a few special cases in which $ch$ is injective.  In the case $d=n-1$
of hypersurfaces, $ch$ is the identity.  In the case $d=1$, in which
$X$ is a point set with multiplicity,
$ch(X)$ is a sum of reflected tropical hyperplanes with multiplicity, 
from which $X$ is easily recoverable.  
Furthermore, Conjecture~\ref{cj:kernel} below would imply 
restrictions on the rays in any one-dimensional tropical fan cycle in $\ker ch$,
and one can check that no cycle with these restrictions lies in $\ker ch$.

Example~\ref{ex:kernel element} provides an explicit tropical fan cycle 
in~$\ker ch$ in the least case, $(d,n)=(3,5)$,
not among those just mentioned.  First we introduce the fan
on which the example depends, which seems to be of critical importance 
to the behaviour of $\ker ch$ in general.  

Let $\mc A_n\subseteq\R^{n-1}$ be the fan in $N_\R$ consisting
of the cones $\R_{\geq0}\{\eN{J_1},\ldots,\eN{J_i}\}$ for all chains of subsets
$$\emptyset\subsetneq J_1\subsetneq\cdots\subsetneq J_i\subsetneq [n].$$
This fan $\mc A_n$ makes many appearances in combinatorics.  
It is the normal fan of the permutahedron, and by Theorem~\ref{r:GGMS}
also the common refinement of all normal fans of matroid polytopes.  
Its face poset is the order poset of the boolean lattice.
Moreover, its codimension 1 skeleton is supported on the union of
the hyperplanes $\{\{x_i=x_j\}:i\neq j\in[n]\}$ of the {\em type A reflection arrangement},
i.e.\ the {\em braid arrangement}.  

As in Section~\ref{ssec:tropical cycles}, the ring $ Z^\fan(\mc A_n)$ is the Chow cohomology ring
of the toric variety associated to $\Sigma$.  This toric variety
is the closure of the torus orbit of a generic point in the complete flag variety
(which, to say it differently, is $\bb P^{n-1}$ blown up along all the
coordinate subspaces).
The cohomology of this variety has been studied by Stembridge~\cite{Stembridge}.
We have that $\dim  Z^\fan(\mc A_n) = n!$, and
$\dim ( Z^\fan)^k(\mc A_n)$ is the \emph{Eulerian number} $E(n,k)$, 
the number of permutations of~$[n]$ with $k$ descents.

For any cone $\sigma=\Rp\{\eN{J_1},\ldots,\eN{J_d}\}$ of~$\mc A_n$, and 
any orthant ${\sigma_{J'}}^\refl = \Rp\{-\eN{j} : j\in J'\}$,
the Minkowski sum $\sigma + {\sigma_{J'}}^\refl$ 
is again a union of cones of~$\mc A_n$.
Therefore $ch( Z^\fan_d(\mc A_n))\subseteq ( Z^\fan)^1(\mc A_n)$ always, and
we find nontrivial elements of $\ker ch$ whenever the dimension
of $ Z^\fan_d(\mc A_n)$ exceeds that of $( Z^\fan)^1(\mc A_n)$, i.e.\ when
$E(n,n-d)>E(n,1)$, equivalently when $2<d<n-1$.  

\begin{example}\label{ex:kernel element}
\begin{figure}
\centering
\includegraphics[width=0.3\textwidth]{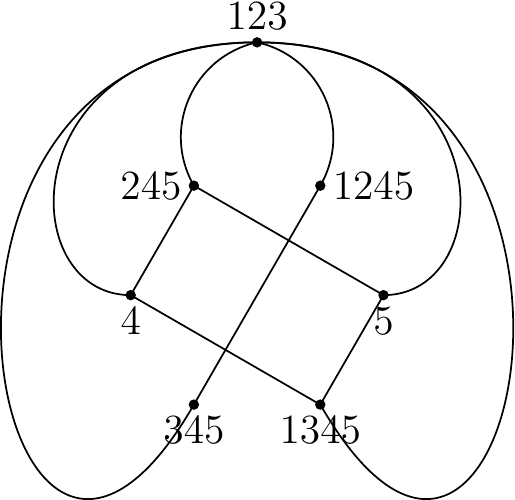}
\hspace{0.5in} 
\includegraphics[width=0.3\textwidth]{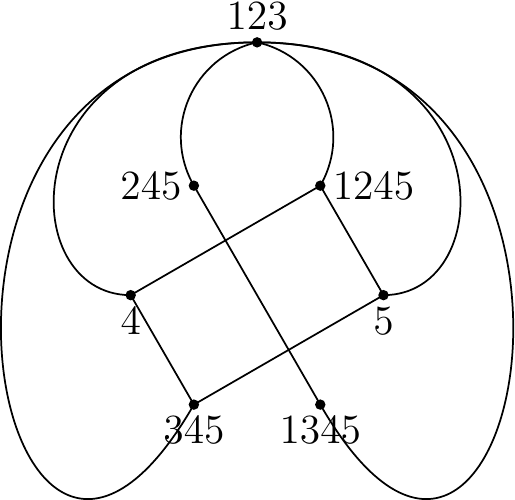}
\caption{Two tropical surfaces with the same Chow hypersurface.
A point labelled $j_1\cdots j_k$ represents the ray
$\R_{\geq0}\eN{\{j_1,\ldots,j_k\}}$.}
\label{fig:kernel element}
\end{figure}

For $(d,n)=(3,5)$, we have $E(5,5-3) = 66 > 26 = E(5,1)$, and the kernel
of $ch$ restricted to $ Z^\fan_2(\mc A_5)$ is 40-dimensional.
Two tropical varieties in $ Z^\fan_2(\mc A_5)$ within $N_\R=\R^4$
with equal Chow hypersurfaces are 
depicted in Figure~\ref{fig:kernel element}. 
As one often does, we have dropped one dimension in the drawing
by actually drawing the intersections of these
2-dimensional tropical fans with a sphere centered at the origin in $\R^4$, 
which are graphs in $\R^3$.  
The difference of these varieties is an actual element of $\ker ch$,
involving the six labelled rays other than 123, which form an octahedron.
\end{example}



The property of $\mc A_n$ that this example exploits appears to be essentially unique:
this is part~(a) of the next conjecture.  
This property, together with experimentation with fan varieties of low degree in low ambient dimension,
also suggests part~(b).

\begin{conjecture}\label{cj:kernel}\mbox{}
\begin{enumerate}
\renewcommand{\labelenumi}{(\alph{enumi})}
\item Let $\Sigma$ be a complete fan such that the stable Minkowski sum of
any cone of $\Sigma$ and any ray $\R_{\geq0}(-\eN{i})$ is a sum of cones of $\Sigma$.
Then $\mc A_n$ is a refinement of~$\Sigma$.
\item The kernel of the restriction of $ch$ to fan varieties 
is generated by elements of~$ Z^\fan(\mc A_n)$. 
\end{enumerate}
\end{conjecture}

\section*{Acknowledgements}
The author thanks Johannes Rau for helpful discussion,
Bernd Sturmfels, Federico Ardila and Eric Katz for close readings
and useful suggestions, and David Speyer for comments.

\end{document}